\documentclass[11pt]{article}
\usepackage{epsfig}

\def\be{\begin{equation}}
\def\ee{\end{equation}}
\def\bea{\begin{eqnarray}}
\def\eea{\end{eqnarray}}
\def\bes{\begin{eqnarray*}}
\def\ees{\end{eqnarray*}}

\def\nn{\nonumber}
\def\<{\langle}
\def\>{\rangle}
\def\lb{\label}
\def\bs{\setminus}
\def\pt{\partial}

\def\R{{\bf R}}
\def\C{{\bf C}}
\def\Z{{\bf Z}}
\def\N{{\bf N}}
\def\U{{\bf U}}
\def\Q{{\bf Q}}

\def\ga{{\gamma}}

\def\ka{{\kappa}}
\def\th{{\theta}}

\def\lm{{\lambda}}
\def\Lm{{\Lambda}}

\def\rank{{\rm rank}}
\def\Sp{{\rm Sp}}
\def\mod{{\rm mod}}

\def\ol#1{\overline{#1}}  

\def\hb{\vrule height0.18cm width0.14cm $\,$}

\def\ol#1{\overline{#1}}  

\title{ On a conjecture of Anosov }
\author{Wei Wang\thanks{Partially supported by National Natural
Science Foundation of China No.10801002, China Postdoctoral Science
Foundation No.200801021, Foundation for the Author of National
Excellent Doctoral Dissertation of PR China No. 201017.
E-mail: alexanderweiwang@yahoo.com.cn, wangwei@math.pku.edu.cn  }\\
School of Mathematical Science \\ Peking University, Beijing 100871 \\
PEOPLES REPUBLIC OF CHINA \\ }
\date{May 1st, 2011 }

\begin{document}

\maketitle

\begin{abstract}
{\it In this paper, we prove that for every bumpy Finsler  $n$-sphere
$(S^n,\,F)$ with reversibility $\lambda$ and flag
curvature $K$ satisfying
$\left(\frac{\lambda}{\lambda+1}\right)^2<K\le 1$,  there
exist $2[\frac{n+1}{2}]$ prime closed geodesics.
This gives a confirmed answer to a conjecture of D. V. Anosov \cite{Ano} in 1974 for a generic case. }
\end{abstract}



{\bf Running head}: On a conjecture of Anosov

\renewcommand{\theequation}{\thesection.\arabic{equation}}
\renewcommand{\thefigure}{\thesection.\arabic{figure}}

\setcounter{equation}{0}
\section{Introduction and main results}

This paper is devoted to a study on closed geodesics on Finsler
spheres.  Let us recall firstly the definition of the Finsler
metrics.

{\bf Definition 1.1.} (cf. \cite{She}) {\it Let $M$ be a finite
dimensional manifold. A function $F:TM\to [0,+\infty)$ is a {\rm
Finsler metric} if it satisfies

(F1) $F$ is $C^{\infty}$ on $TM\bs\{0\}$,

(F2) $F(x,\lm y) = \lm F(x,y)$ for all $y\in T_xM$, $x\in M$, and
$\lm>0$,

(F3) For every $y\in T_xM\bs\{0\}$, the quadratic form
$$ g_{x,y}(u,v) \equiv
         \frac{1}{2}\frac{\pt^2}{\pt s\pt t}F^2(x,y+su+tv)|_{t=s=0},
         \qquad \forall u, v\in T_xM, $$
is positive definite.

In this case, $(M,F)$ is called a {\rm Finsler manifold}. $F$ is
{\rm reversible} if $F(x,-y)=F(x,y)$ holds for all $y\in T_xM$ and
$x\in M$. $F$ is {\rm Riemannian} if $F(x,y)^2=\frac{1}{2}G(x)y\cdot
y$ for some symmetric positive definite matrix function $G(x)\in
GL(T_xM)$ depending on $x\in M$ smoothly. }

A closed curve in a Finsler manifold is a closed geodesic if it is
locally the shortest path connecting any two nearby points on this
curve (cf. \cite{She}). As usual, on any Finsler n-sphere
$S^n=(S^n,F)$, a closed geodesic $c:S^1=\R/\Z\to S^n$ is {\it prime}
if it is not a multiple covering (i.e., iteration) of any other
closed geodesics. Here the $m$-th iteration $c^m$ of $c$ is defined
by $c^m(t)=c(mt)$. The inverse curve $c^{-1}$ of $c$ is defined by
$c^{-1}(t)=c(1-t)$ for $t\in \R$. Note that on a non-symmetric
Finsler manifold, the inverse curve of a closed geodesic is
not a closed geodesic in general. We call two prime closed geodesics
$c$ and $d$ {\it distinct} if there is no $\th\in (0,1)$ such that
$c(t)=d(t+\th)$ for all $t\in\R$. We shall omit the word {\it
distinct} when we talk about more than one prime closed geodesic.
On a symmetric Finsler (or Riemannian) $n$-sphere, two closed geodesics
$c$ and $d$ are called { \it geometrically distinct} if $
c(S^1)\neq d(S^1)$, i.e., their image sets in $S^n$ are distinct.

For a closed geodesic $c$ on $(S^n,\,F)$, denote by $P_c$
the linearized Poincar\'{e} map of $c$. Then $P_c\in \Sp(2n-2)$ is symplectic.
For any $M\in \Sp(2k)$, we define the {\it elliptic height } $e(M)$
of $M$ to be the total algebraic multiplicity of all eigenvalues of
$M$ on the unit circle $\U=\{z\in\C|\; |z|=1\}$ in the complex plane
$\C$. Since $M$ is symplectic, $e(M)$ is even and $0\le e(M)\le 2k$.
A closed geodesic $c$ is called {\it elliptic} if $e(P_c)=2(n-1)$, i.e., all the
eigenvalues of $P_c$ locate on $\U$; {\it hyperbolic} if $e(P_c)=0$, i.e., all the
eigenvalues of $P_c$ locate away from $\U$;
{\it non-degenerate} if $1$
is not an eigenvalue of $P_c$. A Finsler sphere $(S^n,\,F)$
is called {\it bumpy} if all the closed geodesics on it are
non-degenerate.

Following H-B. Rademacher in
\cite{Rad3}, the reversibility $\lambda=\lambda(M,\,F)$ of a compact
Finsler manifold $(M,\,F)$ is defied to be
$$\lambda:=\max\{F(-X)\,|\,X\in TM, \,F(X)=1\}\ge 1.$$

It was quite surprising when Katok \cite{Kat} in 1973 found some non-reversible
Finsler metrics on CROSS with only finitely many prime closed geodesics and all closed
geodesics are non-degenerate and elliptic. The smallest number of closed geodesics
 on $S^n$ that one obtains in these examples is $2[\frac{n+1}{2}]$ (cf. \cite{Zil}).
Then D. V. Anosov in I.C.M. of 1974 conjectured that
the lower bound of the number of closed geodesics on
any Finsler sphere $(S^n,\,F)$ should be $2[\frac{n+1}{2}]$,
i.e., the number of closed geodesics in Katok's example.

We can show that under some conditions this conjecture is true,
i.e., the following main result of the paper.

{\bf Theorem 1.2.} {\it For every bumpy Finsler  $n$-sphere
$(S^n,\,F)$  with reversibility $\lambda$ and flag
curvature $K$ satisfying
$\left(\frac{\lambda}{\lambda+1}\right)^2<K\le 1$,  there
exist $2[\frac{n+1}{2}]$ prime closed geodesics.   }

We can obtain a stability result.

{\bf Theorem 1.3.} {\it For every bumpy Finsler  $n$-sphere
$(S^n,\,F)$  with reversibility $\lambda$ and flag
curvature $K$ satisfying
$\left(\frac{3\lambda}{2(\lambda+1)}\right)^2<K\le 1$,  there
exist two elliptic prime closed geodesics provided
the number of prime closed geodesics on $(S^n,\,F)$ is finite.  }

{\bf Remark 1.4.}
In \cite{BTZ2}, W. Ballmann, G. Thorbergsson and
W. Ziller proved that for a Riemannian metric on $S^n$ with sectional curvature $1/4\le K\le 1$ there exist
$g(n)$ geometrically distinct closed geodesics, and $\frac{n(n+1)}{2}$
geometrically distinct closed geodesics if the metric is bumpy.
In \cite{BaL}, V. Bangert and Y. Long proved that on any Finsler 2-sphere $(S^2, F)$,
there exist at least two prime closed geodesics, which solves
Anosov's conjecture for the $S^2$ case..
In \cite{LoW2} of Y. Long and the author, they further
proved the existence of at least two irrationally elliptic
prime closed geodesics on every Finsler  $2$-sphere $(S^2,\,F)$
provided the number of prime closed geodesics is finite.
In \cite{Rad4}, H.-B. Rademacher studied the existence and
stability of closed geodesics on positively curved Finsler manifolds.
In a series papers \cite{LoD}, \cite{DuL1}-\cite{DuL3} of Y. Long and H. Duan,
they proved there exist two prime closed geodesics on
any compact simply connected Finsler or Riemannian manifold.
In \cite{Rad5}, H.-B. Rademacher proved there exist two prime closed geodesics on
any bumpy $n$-sphere. In \cite{Wang}, the author proved there exist three prime
closed geodesics on any $(S^3, F)$ satisfying
$(\frac{\lambda}{\lambda+1})^2<K\le 1$.

Our proof of these theorems contains mainly three ingredients: the
common index jump theorem of Y. Long, Morse theory and the mean index
identity of H.-B. Rademacher. Fix a Finsler metric $F$
on $S^n$. Let $\Lm=\Lm S^n$ be the free loop space of $S^n$, which
is a Hilbert manifold. For definition and basic properties of
$\Lm$, we refer readers to \cite{Kli2} and \cite{Kli3}. Let
${E}(c)=\frac{1}{2}\int_0^1F(\dot{c}(t))^2dt$ be the energy
functional on $\Lm$. In this paper for $\ka\in \R$ we denote by
\be \Lm^{\ka}=\{d\in \Lm\,|\,E(d)\le \ka\},    \lb{1.1}\ee
and consider the quotient space $\Lm/S^1$. Since the
energy functional ${E}$ is $S^1$-invariant, the negative
gradient flow of $E$ induce a flow on $\Lm/S^1$, so we can apply
Morse theory on $\Lm/S^1$. By a result of H-B. Rademacher in
\cite{Rad1} of 1989, we get the Morse series of the space pair
$(\Lm/S^1, \Lm^0/S^1)$ with rational coefficients. The reason we
use $(\Lm/S^1, \Lm^0/S^1)$ instead of $(\Lm, \Lm^0)$ is that
the Morse series of the first is lacunary.

In this paper, let $\N$, $\N_0$, $\Z$, $\Q$, $\R$, and $\C$ denote
the sets of natural integers, non-negative integers, integers,
rational numbers, real numbers, and complex numbers respectively.
We use only singular homology modules with $\Q$-coefficients.
For terminologies in algebraic topology we refer to \cite{GrH}.
For $k\in\N$, we denote by $\Q^k$ the direct sum $\Q\oplus\cdots\oplus\Q$ of
$k$ copies of $\Q$ and $\Q^0=0$. For an $S^1$-space $X$,
we denote by $\overline{X}$ the quotient space $X/S^1$.
We define the functions
\be \left\{\matrix{[a]=\max\{k\in\Z\,|\,k\le a\}, &
           \mathcal{E}(a)=\min\{k\in\Z\,|\,k\ge a\} , \cr
    \varphi(a)=\mathcal{E}(a)-[a], &\{a\}=a-[a].  \cr}\right. \lb{1.2}\ee
Especially, $\varphi(a)=0$ if $ a\in\Z\,$, and $\varphi(a)=1$ if $
a\notin\Z\,$.

\setcounter{equation}{0}
\section{Variational structures of closed geodesics}

In this section, we review the variational structures of closed geodesic,
all the details can be found in \cite{Rad2} or \cite{BaL}.

On a compact Finsler manifold $(M,F)$, we choose an auxiliary Riemannian
metric. This endows the space $\Lambda=\Lambda M$ of $H^1$-maps
$\gamma:S^1\rightarrow M$ with a natural structure of Riemannian
Hilbert manifolds on which the group $S^1=\R/\Z$ acts continuously
by isometries, cf. \cite{Kli2}, Chapters 1 and 2. This action is
defined by translating the parameter, i.e.
$$ (s\cdot\gamma)(t)=\gamma(t+s) \qquad $$
for all $\gamma\in\Lm$ and $s,t\in S^1$.
The Finsler metric $F$ defines an energy functional $E$ and a length
functional $L$ on $\Lambda$ by
\be E(\gamma)=\frac{1}{2}\int_{S^1}F(\dot{\gamma}(t))^2dt,
 \quad L(\gamma) = \int_{S^1}F(\dot{\gamma}(t))dt.  \lb{2.1}\ee
Both functionals are invariant under the $S^1$-action. The critical points
of $E$ of positive energies are precisely the closed geodesics $c:S^1\to M$
of the Finsler structure. If $c\in\Lambda$ is a closed geodesic then $c$ is
a regular curve, i.e. $\dot{c}(t)\not= 0$ for all $t\in S^1$, and this implies
that the second differential $E''(c)$ of $E$ at $c$ exists.
As usual we define the index $i(c)$ of $c$ as the maximal dimension of
subspaces of $T_c \Lambda$ on which $E^{\prime\prime}(c)$ is negative definite, and the
nullity $\nu(c)$ of $c$ so that $\nu(c)+1$ is the dimension of the null
space of $E^{\prime\prime}(c)$.

For $m\in\N$ we denote the $m$-fold iteration map
$\phi^m:\Lambda\rightarrow\Lambda$ by \be \phi^m(\ga)(t)=\ga(mt)
\qquad \forall\,\ga\in\Lm, t\in S^1. \lb{2.2}\ee We also use the
notation $\phi^m(\gamma)=\gamma^m$. For a closed geodesic $c$, the
mean index is defined to be: \be
\hat{i}(c)=\lim_{m\rightarrow\infty}\frac{i(c^m)}{m}. \lb{2.3}\ee

If $\gamma\in\Lambda$ is not constant then the multiplicity
$m(\gamma)$ of $\gamma$ is the order of the isotropy group $\{s\in
S^1\mid s\cdot\gamma=\gamma\}$. If $m(\gamma)=1$ then $\gamma$ is
called {\it prime}. Hence $m(\gamma)=m$ if and only if there exists a
prime curve $\tilde{\gamma}\in\Lambda$ such that
$\gamma=\tilde{\gamma}^m$.

For a closed geodesic $c$ we set
$$ \Lm(c)=\{\ga\in\Lm\mid E(\ga)<E(c)\}. $$
If $A\subseteq\Lm$ is invariant under some subgroup $\Gamma$ of $S^1$,
we denote by $A/\Gamma$ the
quotient space of $A$ with respect to the action of $\Gamma$.

Using singular homology with rational coefficients we will
consider the following critical $\Q$-module of a closed geodesic
$c\in\Lambda$:
\be \overline{C}_*(E,c)
   = H_*\left((\Lm(c)\cup S^1\cdot c)/S^1,\Lm(c)/S^1\right). \lb{2.4}\ee

We call a closed geodesic satisfying the isolation condition, if
the following holds:

{\bf (Iso)  For all $m\in\N$ the orbit $S^1\cdot c^m$ is an
isolated critical orbit of $E$. }

Note that if the number of prime closed geodesics on a Finsler manifold
is finite, then all the closed geodesics satisfy (Iso).

The following Propositions were proved in \cite{Rad2} and \cite{BaL}.

{\bf Proposition 2.1.} (cf. Satz 6.11 of \cite{Rad2} or Proposition
3.12 of \cite{BaL}) {\it Let $c$ be a prime closed geodesic on a
bumpy Finsler manifold $(M,F)$ satisfying (Iso). Then we have
\bea \overline{C}_q( E,c^m) = \left\{\matrix{
     \Q, &\quad {\it if}\;\; i(c^m)-i(c)\in 2\Z,\;{\it and}\;
                   q=i(c^m)\;  \cr
     0, &\quad {\it otherwise}. \cr}\right.  \lb{2.5}\eea}

Now we briefly describe the relative homological
structure of the quotient space $\overline{\Lm}\equiv
\overline{\Lm} S^n$. Here we have $\ol{\Lm}^0=\ol{\Lambda}^0S^n
=\{{\rm constant\;point\;curves\;in\;}S^n\}\cong S^n$.

{\bf Theorem 2.2.} (H.-B. Rademacher, Theorem 2.4 and Remark 2.5 of \cite{Rad1})
{\it We have the Poincar\'e series

(i) When $n=2k+1$ is odd
\bea
P(\ol{\Lm}S^n,\ol{\Lm}^0S^n)(t)
&=&t^{n-1}\left(\frac{1}{1-t^2}+\frac{t^{n-1}}{1-t^{n-1}}\right)\nn\\
&=& t^{2k}\left(\frac{1}{1-t^2}+\frac{t^{2k}}{1- t^{2k}}\right).
\lb{2.6}\eea
Thus for $q\in\Z$ and $l\in\N_0$, we have
\bea {b}_q &=& {b}_q(\ol{\Lm}S^n,\ol{\Lm}^0 S^n)\nn\\
 &=&\rank H_q(\ol{\Lm} S^n,\ol{\Lm}^0 S^n )\nn\\
 &=& \;\;\left\{\matrix{
    2,&\quad {\it if}\quad q\in \{4k+2l,\quad l=0\;\mod\; k\},  \cr
    1,&\quad {\it if}\quad q\in \{2k\}\cup\{2k+2l,\quad l\neq 0\;\mod\; k\},  \cr
    0 &\quad {\it otherwise}. \cr}\right. \lb{2.7}\eea

(ii) When $n=2k$ is even
\bea
P(\ol{\Lm}S^n,\ol{\Lm}^0S^n)(t)
&=&t^{n-1}\left(\frac{1}{1-t^2}+\frac{t^{n(m+1)-2}}{1-t^{n(m+1)-2}}\right)
\frac{1-t^{nm}}{1-t^n}\nn\\
&=& t^{2k-1}\left(\frac{1}{1-t^2}+\frac{t^{4k-2}}{1- t^{4k-2}}\right),
\lb{2.8}\eea
where $m=1$ by Theorem 2.4 of \cite{Rad1}. Thus for $q\in\Z$ and $l\in\N_0$, we have}
\bea {b}_q &=& {b}_q(\ol{\Lm}S^n,\ol{\Lm}^0 S^n)\nn\\
 &=&\rank H_q(\ol{\Lm} S^n,\ol{\Lm}^0 S^n )\nn\\
 &=& \;\;\left\{\matrix{
    2,&\quad {\it if}\quad q\in \{6k-3+2l,\quad l=0\;\mod\; 2k-1\},  \cr
    1,&\quad {\it if}\quad q\in \{2k-1\}\cup\{2k-1+2l,\quad l\neq 0\;\mod\; 2k-1\},  \cr
    0 &\quad {\it otherwise}. \cr}\right. \lb{2.9}\eea

We have the following version of the Morse inequality.

{\bf Theorem 2.3.} (Theorem 6.1 of \cite{Rad2}) {\it Suppose that there exist
only finitely many prime closed geodesics $\{c_j\}_{1\le j\le p}$ on $(M, F)$,
and $0\le a<b\le \infty$ are regular values of the energy functional $E$.
Define for each $q\in\Z$,
\bea
{M}_q(\ol{\Lm}^b,\ol{\Lm}^a)
&=& \sum_{1\le j\le p,\;a<E(c^m_j)<b}\rank{\ol{C}}_q(E, c^m_j ) \nn\\
{b}_q(\ol{\Lm}^{b},\ol{\Lm}^{a})
&=& \rank H_q(\ol{\Lm}^{b},\ol{\Lm}^{a}). \nn\eea
Then there holds }
\bea
M_q(\ol{\Lm}^{b},\ol{\Lm}^{a}) &-& M_{q-1}(\ol{\Lm}^{b},\ol{\Lm}^{a})
    + \cdots +(-1)^{q}M_0(\ol{\Lm}^{b},\ol{\Lm}^{a}) \nn\\
&\ge& b_q(\ol{\Lm}^{b},\ol{\Lm}^{a}) - b_{q-1}(\ol{\Lm}^{b},\ol{\Lm}^{a})
   + \cdots + (-1)^{q}b_0(\ol{\Lm}^{b},\ol{\Lm}^{a}), \lb{2.10}\\
{M}_q(\ol{\Lm}^{b},\ol{\Lm}^{a}) &\ge& {b}_q(\ol{\Lm}^{b},\ol{\Lm}^{a}).\lb{2.11}
\eea

\setcounter{equation}{0}
\section{Classification of closed geodesics on $S^n$ }

Let $c$ be a closed geodesic on a Finsler n-sphere $S^n=(S^n,\,F)$.
Denote the linearized Poincar\'e map of $c$ by $P_c\in\Sp(2n-2)$.
Then $P_c$ is a symplectic matrix.
Note that the index iteration formulae in \cite{Lon3} of 2000 (cf. Chap. 8 of
\cite{Lon4}) work for Morse indices of iterated closed geodesics (cf.
\cite{LiL}, Chap. 12 of \cite{Lon4}). Since every closed geodesic
on a  sphere must be orientable. Then by Theorem 1.1 of \cite{Liu}
of C. Liu (cf. also \cite{Wil}), the initial Morse index of a closed geodesic
$c$ on a $n$-dimensional Finsler  sphere coincides with the index of a
corresponding symplectic path introduced by C. Conley, E. Zehnder, and Y. Long
in 1984-1990 (cf. \cite{Lon4}).

As in \S 1.8 of \cite{Lon4}, define the homotopy component $\Omega^0 (P_c)$
of $P_c$ to be the path component of $\Omega (P_c)$, where
\bea \Omega (P_c) =\{N\in
Sp(2n-2) \mid &&\sigma (N)\cap U=\sigma (P_c)\cap U, \;and\nn\\
&&\nu_\lambda(N)=\nu_\lambda(P_c) \;\forall \lambda\in \sigma
(P_c)\cap U\}.\lb{3.1}\eea
The next theorem is due to Y. Long (cf. Theorem 8.3.1 and Corollary 8.3.2 of \cite{Lon4}).

{\bf Theorem 3.1.} {\it Let $\gamma\in\{\xi\in
C([0,\tau], Sp(2n))\mid \xi(0)=I\}$, Then there exists a path $f\in
C([0,1],\Omega^0(\gamma(\tau))$ such that $f(0)=\gamma(\tau)$ and
\bea f(1)=&&N_1(1,1)^{\diamond p_-} \diamond I_{2p_0}\diamond
N_1(1,-1)^{\diamond p_+}
\diamond N_1(-1,1)^{\diamond q_-} \diamond (-I_{2q_0})\diamond
N_1(-1,-1)^{\diamond q_+}\nn\\
&&\diamond R(\theta_1)\diamond\cdots\diamond R(\theta_r)
\diamond N_2(\omega_1, u_1)\diamond\cdots\diamond N_2(\omega_{r_*}, u_{r_*}) \nn\\
&&\diamond N_2(\lm_1, v_1)\diamond\cdots\diamond N_2(\lm_{r_0}, v_{r_0})
\diamond M_0 \lb{3.2}\eea
where $ N_2(\omega_j, u_j) $s are
non-trivial and   $ N_2(\lm_j, v_j)$s  are trivial basic normal
forms; $\sigma (M_0)\cap U=\emptyset$; $p_-$, $p_0$, $p_+$, $q_-$,
$q_0$, $q_+$, $r$, $r_*$ and $r_0$ are non-negative integers;
$\omega_j=e^{\sqrt{-1}\alpha_j}$, $
\lambda_j=e^{\sqrt{-1}\beta_j}$; $\theta_j$, $\alpha_j$, $\beta_j$
$\in (0, \pi)\cup (\pi, 2\pi)$; these integers and real numbers
are uniquely determined by $\gamma(\tau)$. Then using the
functions defined in (\ref{1.2}).
\bea i(\gamma, m)=&&m(i(\gamma,
1)+p_-+p_0-r)+2\sum_{j=1}^r\mathcal{E}\left(\frac{m\theta_j}{2\pi}\right)-r
-p_--p_0\nn\\&&-\frac{1+(-1)^m}{2}(q_0+q_+)+2\left(
\sum_{j=1}^{r_*}\varphi\left(\frac{m\alpha_j}{2\pi}\right)-r_*\right),
\lb{3.3}\eea
\bea \nu(\gamma, m)=&&\nu(\gamma,
1)+\frac{1+(-1)^m}{2}(q_-+2q_0+q_+)+2(r+r_*+r_0)\nn\\
&&-2\left(\sum_{j=1}^{r}\varphi\left(\frac{m\theta_j}{2\pi}\right)+
\sum_{j=1}^{r_*}\varphi\left(\frac{m\alpha_j}{2\pi}\right)
+\sum_{j=1}^{r_0}\varphi\left(\frac{m\beta_j}{2\pi}\right)\right),\lb{3.4}\eea
\bea \hat i(\gamma, 1)=i(\gamma, 1)+p_-+p_0-r+\sum_{j=1}^r
\frac{\theta_j}{\pi}.\lb{3.5}\eea
Where $N_1(1, \pm 1)=
\left(\matrix{ 1 &\pm 1\cr 0 & 1\cr}\right)$, $N_1(-1, \pm 1)=
\left(\matrix{ -1 &\pm 1\cr 0 & -1\cr}\right)$,
$R(\theta)=\left(\matrix{\cos\th &
                  -\sin\th\cr\sin\th & \cos\th\cr}\right)$,
$ N_2(\omega, b)=\left(\matrix{R(\th) & b
                  \cr 0 & R(\th)\cr}\right)$    with some
$\th\in (0,\pi)\cup (\pi,2\pi)$ and $b=
\left(\matrix{ b_1 &b_2\cr b_3 & b_4\cr}\right)\in\R^{2\times2}$,
such that $(b_2-b_3)\sin\theta>0$, if $ N_2(\omega, b)$ is
trivial; $(b_2-b_3)\sin\theta<0$, if $ N_2(\omega, b)$ is
non-trivial. We have $i(\gamma, 1)$ is odd if $f(1)=N_1(1, 1)$, $I_2$,
$N_1(-1, 1)$, $-I_2$, $N_1(-1, -1)$ and $R(\theta)$; $i(\gamma, 1)$ is
even if $f(1)=N_1(1, -1)$ and $ N_2(\omega, b)$; $i(\gamma, 1)$ can be any
integer if $\sigma (f(1)) \cap \U=\emptyset$.}

\setcounter{equation}{0}
\section{A mean index equality on $(S^n, F)$}

In this section, we recall the mean index equality obtained in
\cite{Rad1}. Suppose that there are only finitely many prime closed
geodesics $\{c_j\}_{1\le j\le p}$ on a bumpy $(S^n, F)$ with $\hat i(c_j)>0$
for $1\le j\le p$.

{\bf Lemma 4.1.} {\it  Let $c$ be a prime closed geodesic on
a bumpy $(S^n,F)$.
Then we have
\bea i(c^{p+2})-i(c^p)\in 2\Z,
     \qquad\forall p\in \N.  \lb{4.1}\eea
 }
{\bf Proof.} This follows directly from Theorem 3.1. \hfill\hb

{\bf Definition 4.2.} {\it Suppose $c$ is a closed
geodesic on $(S^n, F)$. The Euler characteristic $\chi(c^m)$
of $c^m$ is defined by
\bea \chi(c^m)
&\equiv& \chi\left((\Lm(c^m)\cup S^1\cdot c^m)/S^1, \Lm(c^m)/S^1\right), \nn\\
&\equiv& \sum_{q=0}^{\infty}(-1)^q\dim \overline{C}_q( E,c^m).
\lb{4.2}\eea
Here $\chi(A, B)$ denotes the usual Euler characteristic of the space pair $(A, B)$.

The average Euler characteristic $\hat\chi(c)$ of $c$ is defined by }
\be \hat{\chi}(c)=\lim_{N\to\infty}\frac{1}{N}\sum_{1\le m\le N}\chi(c^m).
\lb{4.3}\ee

The following remark shows that $\hat\chi(c)$ is well-defined and is a
rational number.

{\bf Remark 4.3.} By  Proposition 2.1 and Lemma 4.1 we have
\bea
\hat\chi(c) &=&\lim_{N\rightarrow\infty}\frac{1}{N}
  \sum_{1\le m\le N}(-1)^{i(c^m)}\dim \overline{C}_{i(c^m)}( E, c^m)\nn\\
&=&\lim_{s\rightarrow\infty}\frac{1}{2s}
  \sum_{1\le m\le 2\atop  0\le p< s}
  (-1)^{i(c^{2p+m})}\dim \overline{C}_{i(c^{2p+m})}( E, c^m)\nn\\
&=&\frac{1}{2}
  \sum_{1\le m\le 2}
  (-1)^{i(c^{m})}\dim \overline{C}_{i(c^m)}( E, c^m)
=\frac{1}{2}\sum_{1\le m\le 2}\chi(c^m).\lb{4.4}\eea
Therefore $\hat\chi(c)$ is well defined and is a rational number.

The following is the mean index equality of H.-B. Rademacher (Theorem 7.9
in \cite{Rad2}).

{\bf Theorem 4.4.} {\it Suppose that there exist only finitely many prime
closed geodesics $\{c_j\}_{1\le j\le p}$ with $\hat i(c_j)>0$
for $1\le j\le p$ on $(S^n,F)$.
Then the following identity holds
\be  \sum_{1\le j\le p}\frac{\hat\chi(c_j)}{\hat{i}(c_j)}=B(n,1)
=\left\{\matrix{
     \frac{n+1}{2(n-1)}, &\quad n\quad odd,  \cr
     \frac{-n}{2(n-1)}, &\quad n\quad even.\cr}\right.\lb{4.5}\ee
}

\setcounter{equation}{0}
\section{ Proof of the main theorems }

In this section, we give the proofs of Theorems 1.2 and 1.3 by using
the mean index identity in Theorem 4.4, Morse inequality and the index iteration theory developed by
Y. Long and his coworkers.

In the following for the notation introduced in Section 3 we use
specially $M_j=M_j(\ol{\Lm} S^n,\ol{\Lm}^0 S^n)$ and
$b_j=b_j(\ol{\Lm} S^n,\ol{\Lm}^0 S^n)$ for $j=0,1,2,\ldots$.

First note that if the flag curvature $K$ of $(S^n, F)$ satisfies
$\left(\frac{\lambda}{\lambda+1}\right)^2<K\le 1$,
then every nonconstant closed geodesic must satisfy
\bea i(c)&\ge& n-1, \lb{5.1}\\
\hat i(c)&>& n-1, \lb{5.2}\eea
where (\ref{5.1}) follows from Theorem 3 and Lemma 3 of \cite{Rad3},
(\ref{5.2}) follows from Lemma 2 of \cite{Rad4}.
Now it follows from Theorem 2.2 of \cite{LoZ}
(Theorem 10.2.3 of \cite{Lon4}) that
\bea i(c^{m+1})-i(c^m)\ge i(c)-\frac{e(P_c)}{2}\ge 0,\quad\forall m\in\N.\lb{5.3}\eea
Here the last inequality holds by (\ref{5.1}) and the fact that $e(P_c)\le 2(n-1)$.

In the rest of this paper, we will assume the following

{\bf (F) There are only finitely many prime closed geodesics
$\{c_j\}_{1\le j\le p}$ on $(S^n,\,F)$. }

By (\ref{5.2}), we can use the common index jump theorem (Theorem 4.3 of
\cite{LoZ}, Theorem 11.2.1 of \cite{Lon4}) to obtain some
$(N, m_1,\ldots,m_p)\in\N^{p+1}$ such that
\bea
i(c_j^{2m_j}) &\ge& 2N-\frac{e(P_{c_j})}{2}\ge 2N-(n-1), \lb{5.4}\\
i(c_j^{2m_j}) &\le& 2N+\frac{e(P_{c_j})}{2}\le 2N+(n-1), \lb{5.5}\\
i(c_j^{2m_j-m}) &\le& 2N-(i(c_j)+2S^+_{P_{c_j}}(1)-\nu(c_j)),\quad \forall m\in\N. \lb{5.6}\\
i(c_j^{2m_j+m}) &\ge& 2N+i(c_j),\quad\forall m\in\N, \lb{5.7}
\eea
where $S^+_{P_{c_j}}(1)$ denote the splitting number of $c_j$ at $1$.
Since the metric is bumpy, $1$ is not an eigenvalue of $P_{c_j}$
for $1\le j\le p$. Thus we have $S^+_{P_{c_j}}(1)=0$ and $\nu(c_j)=0$.
Hence (\ref{5.6}) becomes
\bea i(c_j^{2m_j-m}) \le 2N-i(c_j),\quad \forall m\in\N.
\lb{5.8}
\eea
Moreover we have
\bea \min\left\{\left\{\frac{m_j\theta}{\pi}\right\},\,1-\left\{\frac{m_j\theta}{\pi}\right\}\right\}
<\delta,\lb{5.9}\eea
whenever $e^{\sqrt{-1}\theta}\in\sigma(P_{c_j})$
and $\delta$ can be chosen as small as we want. More precisely, by Theorem 4.1 of
\cite{LoZ} (in (11.1.10) in Theorem 11.1.1 of \cite{Lon4}, with $D_j=\hat i(c_j)$,
we have
\bea m_j=\left(\left[\frac{N}{M\hat i(c_j)}\right]+\xi_j\right)M,\quad 1\le j\le p,\lb{5.10}\eea
where $\xi_j=0$ or $1$ for $1\le j\le p$ and $M\in \N$.
By (11.1.20) in Theorem 11.1.1 of \cite{Lon4},
for any $\epsilon>0$, we can choose $N$ and $\{\xi_j\}_{1\le j\le p}$ such that
\bea \left|\frac{N}{M\hat i(c_j)}-\left[\frac{N}{M\hat i(c_j)}\right]-\xi_j\right|<\epsilon
<\frac{1}{1+\sum_{1\le j\le p}4M|\hat\chi(c_j)|},\quad 1\le j\le p.\lb{5.11}\eea

{\bf Lemma 5.1.} {\it There  exists a prime
closed geodesic $c_{j_0}$ such that $i(c^{2m_{j_0}}_{j_0})=2N+(n-1)$.
In particular, we have $\ol{C}_{2N+n-1}(E, c^{2m_{j_0}}_{j_0})\neq 0$. }

{\bf Proof.} Suppose the contrary. Then by (\ref{5.5}), we have
\bea i(c_j^{2m_j})< 2N+(n-1), \qquad 1\le j\le p.\lb{5.12}
\eea
Now by (\ref{5.1}), (\ref{5.3}), (\ref{5.7}) and (\ref{5.12}), we have
\bea
i(c_j^{m}) &\le& i(c_j^{2m_j}),\quad\forall m<2m_j, \lb{5.13}\\
i(c_j^{2m_j}) &\le& 2N+n-2, \lb{5.14}\\
i(c_j^{m}) &\ge& 2N+n-1,\quad\forall m>2m_j. \lb{5.15}\eea
By Theorem 4.4, we have
\be  \sum_{1\le j\le p}\frac{\hat\chi(c_j)}{\hat{i}(c_j)}=B(n,1)\in\Q.\lb{5.16}\ee
Note by the proof of Theorem 4.1 of \cite{LoZ} (Theorem 11.1.1 of \cite{Lon4}),
we can require that $N\in\N$  further satisfies (cf. (11.1.22) in \cite{Lon4})
\bea 2NB(n,1)\in\Z.\lb{5.17}\eea
Multiplying both sides of (\ref{5.16}) by $2N$ yields
\be  \sum_{1\le j\le p}\frac{2N\hat\chi(c_j)}{\hat{i}(c_j)}=2NB(n,1).\lb{5.18}\ee

{\bf Claim 1.} {\it We have}
\be  \sum_{1\le j\le p}2m_j\hat\chi(c_j)=2NB(n,1).\lb{5.19}\ee
In fact, by (\ref{5.10}) and (\ref{5.18}), we have
\bea &&2NB(n,1)\nn\\
=&&\sum_{1\le j\le p}\frac{2N\hat\chi(c_j)}{\hat{i}(c_j)}\nn\\
=&&\sum_{1\le j\le p}2\hat\chi(c_j)\left(\left[\frac{N}{M\hat i(c_j)}\right]+\xi_j\right)M
+\sum_{1\le j\le p}2\hat\chi(c_j)\left(\frac{N}{M\hat i(c_j)}-\left[\frac{N}{M\hat i(c_j)}\right]-\xi_j\right)M\nn\\
=&& \sum_{1\le j\le p}2m_j\hat\chi(c_j)+\sum_{1\le j\le p}2M\hat\chi(c_j)\epsilon_j.\lb{5.20}
\eea
By  (\ref{4.4}) we have
\bea 2m_j\hat\chi(c_j)\in\Z,\quad1\le j\le p.\lb{5.21}\eea
Now Claim 1 follows by (\ref{5.11}), (\ref{5.17}), (\ref{5.20}) and  (\ref{5.21}).

{\bf Claim 2.} {\it  We have
\be \sum_{1\le j\le p}2m_j\hat\chi(c_j)=M_0-M_1+M_2-\cdots+(-1)^{2N+n-2}M_{2N+n-2},\lb{5.22}\ee
where  $M_q\equiv M_q(\ol{\Lm},\ol{\Lm}^0)$ for $q\in\Z$.
}

In fact, by definition, the right hand side of (\ref{5.22}) is
\bea RHS=\sum_{q\le 2N+n-2\atop m\ge 1,\; 1\le j\le p}(-1)^q\dim\ol{C}_q(E, c^m_j).\lb{5.23}
\eea
By (\ref{5.13})-(\ref{5.15}) and Proposition 2.1, we have
\bea RHS&=&\sum_{1\le j\le p,\;1\le m\le 2m_j\atop q\le 2N+n-2}(-1)^q\dim\ol{C}_q(E, c^m_j),\lb{5.24}\\
&=&\sum_{1\le j\le p,\;1\le m\le 2m_j}\chi(c_j^m),\lb{5.25}
\eea
where the second equality follows from (\ref{4.2}).

By  (\ref{4.4}), we have
\bea  \sum_{1\le m\le 2m_j}\chi(c_j^m)
&=&\sum_{0\le s< m_j,\,1\le m\le 2}\chi(c_j^{2s+m})\nn\\
&=&m_j\sum_{1\le m\le 2}\chi(c_j^m)\nn\\
&=& 2m_j\hat\chi(c_j),\lb{5.26}\eea
This proves Claim 2.

Now we  consider the following two
cases according to the parity of $n$.

{\bf Case 1. } $n=2k+1$ is odd.

In this case, we have by (\ref{4.5})
\be B(n,1)=\frac{n+1}{2(n-1)}=\frac{k+1}{2k}.\lb{5.27}\ee
By the proof of Theorem 4.1 of \cite{LoZ} (Theorem 11.1.1 of \cite{Lon4}),
we may  further assume  $N=2ks$ for some $s\in\N$.

Thus by (\ref{5.19}), (\ref{5.22}) and (\ref{5.27}), we have
\be M_0-M_1+M_2-\cdots+(-1)^{2N+n-2}M_{2N+n-2}=2s(k+1).\lb{5.28}
\ee
On the other hand, we have by (\ref{2.7})
\bea &&b_0-b_1+b_2-\cdots+(-1)^{2N+n-2}b_{2N+n-2}\nn\\
=&&b_{2k}+(b_{2k+2}+\cdots+b_{4k}+\cdots+b_{4sk+2}+\cdots+b_{4sk+2k})-b_{4sk+2k}\nn\\
=&& 1+2s(k-1+2)-2\nn\\
=&& 2s(k+1)-1.\lb{5.29}
\eea
In fact, we cut off the sequence $\{b_{2k+2},\ldots,b_{4sk+2k}\}$ into
$2s$ pieces, each of them contains  $k$ terms. Moreover, each piece contain
$1$ for $k-1$ times and  $2$ for one time. Thus (\ref{5.29}) holds.

Now by (\ref{5.28}), (\ref{5.29}) and Theorem 2.3, we have
\bea -2s(k+1)&=&M_{2N+n-2}-M_{2N+n-3}+\cdots+M_1-M_0\nn\\
&\ge&b_{2N+n-2}-b_{2N+n-3}+\cdots+b_1-b_0\nn\\
&=&-(2s(k+1)-1).\lb{5.30}
\eea
This contradiction yields the lemma for $n$ being odd.

{\bf Case 2. } $n=2k$ is even.

In this case, we have by (\ref{4.5})
\be B(n,1)=\frac{-n}{2(n-1)}=\frac{-k}{2k-1}.\lb{5.31}\ee
By the proof of Theorem 4.1 of \cite{LoZ} (Theorem 11.1.1 of \cite{Lon4}),
we may  further assume $N=(2k-1)s$ for some $m\in\N$.

Thus by (\ref{5.19}), (\ref{5.22}) and (\ref{5.31}), we have
\be M_0-M_1+M_2-\cdots+(-1)^{2N+n-2}M_{2N+n-2}=-2sk.\lb{5.32}
\ee
On the other hand, we have by (\ref{2.9})
\bea &&b_0-b_1+b_2-\cdots+(-1)^{2N+n-2}b_{2N+n-2}\nn\\
=&&-b_{2k-1}-(b_{2k+1}+\cdots+b_{6k-3}+\cdots+b_{(s-1)(4k-2)+2k+1}+\cdots+b_{s(4k-2)+2k-1})\nn\\
&&+b_{s(4k-2)+2k-1}\nn\\
=&& -1-s(2k-2+2)+2\nn\\
=&& -2sk+1.\lb{5.33}
\eea
In fact, we cut off the sequence $\{b_{2k+1},\ldots,b_{s(4k-2)+2k-1}\}$ into
$s$ pieces, each of them contains  $2k-1$ terms. Moreover, each piece contain
$1$ for $2k-2$ times and  $2$ for one time. Thus (\ref{5.33}) holds.

Now by (\ref{5.32}), (\ref{5.33}) and Theorem 2.3, we have
\bea -2sk&=&M_{2N+n-2}-M_{2N+n-3}+\cdots+M_1-M_0\nn\\
&\ge&b_{2N+n-2}-b_{2N+n-3}+\cdots+b_1-b_0\nn\\
&=&-2sk+1.\lb{5.34}
\eea
This contradiction yields the lemma for $n$ being even.
 \hfill\hb

{\bf Lemma 5.2} {\it We have
\bea i(c_j^{2m_j-2})< 2N-(n-1),\lb{5.35}\eea
for $1\le j\le p$.}

{\bf Proof.}  By (\ref{5.3}) and (\ref{5.8}), if $i(c_j)>n-1$,
then (\ref{5.35}) holds. Thus it remains to consider the case
$i(c_j)=n-1$. By (\ref{3.5}) and (\ref{5.2}), we have
\bea \hat i(c_j)&=&i(c_j)+p_-+p_0-r+\sum_{i=1}^r
\frac{\theta_i}{\pi}\nn\\
&=&i(c_j)-r+\sum_{i=1}^r\frac{\theta_i}{\pi}
>n-1,\lb{5.36}\eea
where the second equality follows from $p_-=0=p_0$,
which holds since $c_j$ is non-degenerate. Plugging
$i(c_j)=n-1$ into (\ref{5.36}) yields
\bea  \sum_{i=1}^r\left(\frac{\theta_i}{\pi}-1\right)>0.\lb{5.37}\eea
Hence we can write
\bea P_{c_j}=R(\theta)\diamond M,
\lb{5.38}\eea
for some $\theta\in (\pi, 2\pi)$ and $M\in Sp(2n-4)$.
Thus by Theorem 3.1 and the assumption that $c_j$ is non-degenerate, we have
\bea i(c_j^m)&=&m(i(c_j)-r)+2\sum_{i=1}^r\mathcal{E}\left(\frac{m\theta_i}{2\pi}\right)-r
\nn\\
&=& 2\mathcal{E}\left(\frac{m\theta}{2\pi}\right)-1+i(\gamma, m),
\lb{5.39}\eea
where $\gamma\in\{\xi\in C([0,\tau], Sp(2n-4))\mid \xi(0)=I\}$
satisfies $\gamma(\tau)=M$ and $i(\gamma, 1)=n-2$.
The second equality  follows form the Symplectic additivity of the
index function, cf. Theorem 6.2.7 of \cite{Lon4}.
Note that it follows from Theorem 2.2 of \cite{LoZ}
(Theorem 10.2.3 of \cite{Lon4}) that
\bea i(\gamma, m+1)-i(\gamma , m)\ge i(\gamma, 1)-\frac{e(M)}{2}\ge 0,\quad\forall m\in\N.\lb{5.40}\eea
Here the last inequality holds from $i(\gamma, 1)=n-2$
 and the fact that $e(M)\le 2(n-2)$.
By (\ref{5.3}) and (\ref{5.8}), in order to prove (\ref{5.35}), it is sufficient
to prove
\bea i(c_j^{2m_j-2})<i(c_j^{2m_j-1}).\lb{5.41}\eea
By (\ref{5.39}) and (\ref{5.40}), in order to prove (\ref{5.41}), it is sufficient
to prove
\bea \mathcal{E}\left(\frac{(2m_j-2)\theta}{2\pi}\right)
<\mathcal{E}\left(\frac{(2m_j-1)\theta}{2\pi}\right).
\lb{5.42}\eea
In order to satisfy (\ref{5.42}), it is sufficient
to choose
\bea \delta<\min\left\{\frac{\theta}{\pi}-1,\,1-\frac{\theta}{2\pi}\right\},
\lb{5.43}\eea
where $\theta$ is given by (\ref{5.9}). This proves the lemma
\hfill\hb

{\bf Proof of Theorem 1.2.} Note that by Theorem 2.3, we have
\bea
M_q\equiv M_q(\ol{\Lm},\ol{\Lm}^0)
= \sum_{ m\ge 1,\,1\le j\le p}\rank{\ol{C}}_q(E, c^m_j ).\lb{5.44}\eea

The proof contains three steps:

{\bf Step 1.} {\it We have $p\ge 2[\frac{n+1}{2}]-2$.}

We consider the following two
cases according to the parity of $n$.

{\bf Case 1.1. } $n=2k+1$ is odd.

In this case, as in Case 1 of Lemma 5.1,
we may assume  $N=2ks$ for some $s\in\N$.
Then by Theorem 2.2, we have
\bea b_{2N-(n-1)+2m}=1, \quad 1\le m<k, \; k<m\le n-2; \qquad b_{2N}=2.
\lb{5.45}\eea
Thus by Theorem 2.3, we have
\bea &&M_{2N-(n-1)+2m}\ge b_{2N-(n-1)+2m}=1,\quad 1\le m<k, \; k<m\le n-2;\lb{5.46}\\
&&M_{2N}\ge b_{2N}=2.\lb{5.47}\eea
By Proposition 2.1 and (\ref{5.1}), (\ref{5.7}) and (\ref{5.8}), we have
\bea
M_{2N-(n-1)+2m}=&& \sum_{1\le j\le p}\rank{\ol{C}}_{2N-(n-1)+2m}(E, c^{2m_j}_j )\nn\\
=&&^\#\{j | i(c_j^{2m_j})-i(c_j)\in2\Z,\; i(c_j^{2m_j})=2N-(n-1)+2m\},\lb{5.48}\eea
for $1\le m\le n-2$. Hence we have $p\ge n-1$ by (\ref{5.45})-(\ref{5.48}) and Proposition 2.1.
In fact, only the $2m_j$-th iteration $c_j^{2m_j}$ of $c_j$
contribute at most $1$ to
$$\sum_{1\le m\le n-2}M_{2N-(n-1)+2m}\ge\sum_{1\le m\le n-2}b_{2N-(n-1)+2m} =n-3+2=n-1$$
for each $1\le j\le p$. This yields Step 1 for $n$ being odd.

{\bf Case 1.2. } $n=2k$ is even.

In this case, as in Case 2 of Lemma 5.1, we may assume $N=(2k-1)s$ for some $s\in\N$.
Then by Theorem 2.2, we have
\bea b_{2N-(n-1)+2m}=1, \quad 1\le m \le n-2.
\lb{5.49}\eea
Thus by Theorem 2.3, we have
\bea M_{2N-(n-1)+2m}\ge b_{2N-(n-1)+2m}=1,\quad 1\le m \le n-2.
\lb{5.50}\eea
By Proposition 2.1 and (\ref{5.1}), (\ref{5.7}) and (\ref{5.8}), we have
\bea
M_{2N-(n-1)+2m}=&& \sum_{1\le j\le p}\rank{\ol{C}}_{2N-(n-1)+2m}(E, c^{2m_j}_j )\nn\\
=&&^\#\{j | i(c_j^{2m_j})-i(c_j)\in2\Z,\; i(c_j^{2m_j})=2N-(n-1)+2m\},\lb{5.51}\eea
for $1\le m\le n-2$. Hence as in Case 1, we have $p\ge n-2$ by (\ref{5.49})-(\ref{5.51}) and Proposition 2.1.
This yields Step 1 for $n$ being even.

{\bf Step 2.} {\it  We have $p\ge 2[\frac{n+1}{2}]-1$.}

Denote the $2[\frac{n+1}{2}]-2$ prime closed geodesics obtained in Step 1
by $\{c_{j_1},\ldots,c_{j_{2\left[\frac{n+1}{2}\right]-2}}\}$.
Then by (\ref{5.48}), (\ref{5.51}) and Proposition 2.1,
we have
\bea i(c_{j_k}^{2m_{j_k}})=2N-(n-1)+2\tau_{j_k},
\quad \ol{C}_{2N-(n-1)+2\tau_{j_k}}(E, c^{2m_{j_k}}_{j_k})\neq 0,
\lb{5.52}\eea
for some $1\le \tau_{j_k}\le n-2$ and
$1\le k\le 2[\frac{n+1}{2}]-2$.
On the other hand, by Theorem 5.1, there exists a
closed geodesic $c_{j_0}$ such that
\bea i(c^{2m_{j_0}}_{j_0})=2N+(n-1),\qquad \ol{C}_{2N+n-1}(E, c^{2m_{j_0}}_{j_0})\neq 0.\lb{5.53}\eea
Hence we have $c_{j_0}\notin \{c_{j_1},\ldots,c_{j_{2\left[\frac{n+1}{2}\right]-2}}\}$
by (\ref{5.52}) and (\ref{5.53}). This yields Step 2.

{\bf Step 3.} {\it  We have $p\ge 2[\frac{n+1}{2}]$.}

Denote the $2[\frac{n+1}{2}]-1$ prime closed geodesics obtained in Steps 1 and 2
by $\{c_{j_0},\ldots,c_{j_{2\left[\frac{n+1}{2}\right]-2}}\}$.

By Theorems 2.2 and 2.3,  we have
\bea
M_{n-1}=\sum_{1\le j\le p,\,m\ge 1}\rank{\ol{C}}_{n-1}(E, c^{m}_j )
\ge b_{n-1}=1.\lb{5.54}\eea
Thus it follows from (\ref{5.1}) and (\ref{5.3}) that
there exist at least one closed geodesic $c_j$ such that
$i(c_j)=n-1$. We have the following two cases:

{\bf Case 3.1.} {\it We have $^\#\{j | i(c_j)=n-1\}=1$, i.e.,
there is only one prime closed geodesic which has index $n-1$. }

 Denote the  prime closed geodesic which has index $n-1$
 by $c_\ast$. Then we have
 \bea i(c_{l})>n-1,\qquad l\in\{1,\ldots,p\}\setminus\{\ast\}.
 \lb{5.55}\eea
 Thus it follows from (\ref{5.8}) that
 \bea i(c_l^{2m_l-m}) < 2N-(n-1),\quad \forall m\in\N,\; l\in\{1,\ldots,p\}\setminus\{\ast\}
\lb{5.56}
\eea
By Lemma 5.2 and (\ref{5.3}), we have
 \bea i(c_\ast^{2m_\ast-m}) < 2N-(n-1),\quad \forall m\ge 2.
\lb{5.57}
\eea
Then by (\ref{5.1}), (\ref{5.7}), (\ref{5.52}), (\ref{5.53}),
(\ref{5.56}), (\ref{5.57}) and Proposition 2.1, we have
\bea
\sum_{0\le k\le 2\left[\frac{n+1}{2}\right]-2,\,m\ge 1}\rank{\ol{C}}_{2N-(n-1)}(E, c^{m}_{j_k} )
\le 1.\lb{5.58}\eea
 In fact, the only possible non-zero term is
 $\rank{\ol{C}}_{2N-(n-1)}(E, c_{\ast}^{2m_\ast-1})$.

 By Theorems 2.2 and 2.3,  we have
\bea
M_{2N-(n-1)}=\sum_{1\le j\le p,\,m\ge 1}\rank{\ol{C}}_{2N-(n-1)}(E, c^{m}_j )
\ge b_{2N-(n-1)}=2.\lb{5.59}\eea
 Hence there must be another closed geodesic
 $c_\star\notin \{c_{j_0},\ldots,c_{j_{2\left[\frac{n+1}{2}\right]-2}}\}$
by (\ref{5.58}) and (\ref{5.59}). Especially, we have
\bea i(c_\star^{2m_\star})=2N-(n-1)\lb{5.60}\eea
by (\ref{5.1}), (\ref{5.7}), (\ref{5.56}) and Proposition 2.1.
This yields Case 3.1.

 {\bf Case 3.2.} {\it We have $^\#\{j | i(c_j)=n-1\}> 1$. }

By Proposition 2.1 we have
\bea
M_{n-1}=\sum_{1\le j\le p,\,m\ge 1}\rank{\ol{C}}_{n-1}(E, c^{m}_j )
\ge 2\lb{5.61}\eea
 By (\ref{5.1}), Proposition 2.1, Theorems 2.2 and 2.3 we have
 \bea M_n-M_{n-1}&=&M_n-M_{n-1}+\dots+(-1)^nM_0\nn\\
 &\ge& b_n-b_{n-1}+\dots+(-1)^nb_0= b_n-b_{n-1}=-1.
 \lb{6.62}\eea
 Thus we have
 \bea
M_n=\sum_{1\le j\le p,\,m\ge 1}\rank{\ol{C}}_n(E, c^{m}_j )
\ge 1.\lb{5.63}\eea
Then we must have a prime closed geodesic $c_\star$
with
\bea i(c_\star)=n. \lb{5.64}\eea
In fact, by (\ref{5.63}), we have
\bea\ol{C}_n(E, c^{m}_\star)\neq 0.\lb{5.65}\eea
for some $1\le \star\le p$ and $m\in\N$.
By (\ref{5.1}) and (\ref{5.3}), if  $ i(c_\star)\neq n$,
we must have $i(c_\star)=n-1$. Thus it follows from Proposition 2.1
that $\ol{C}_n(E, c^{m}_\star)=0$ for any $m\in\N$. This
contradict to (\ref{5.65}) and yields (\ref{5.64}).
Now it follows from Proposition 2.1 and (\ref{5.64}) that
\bea \ol{C}_{2N-(n-1)+2m}(E, c^{m}_\star)=0,\qquad \forall m\in\Z.
\lb{5.66}\eea
Thus we have  $c_\star\notin \{c_{j_0},\ldots,c_{j_{2\left[\frac{n+1}{2}\right]-2}}\}$
by (\ref{5.52}) and (\ref{5.53}). This yields Case 3.2.

The proof of Theorem 1.2 is complete. \hfill\hb

{\bf Proof of Theorem 1.3.}
We have two cases according to Step 3 of the above proof.

{\bf Case 1.} { \it If Case 3.1 holds.}

In this case, the closed geodesic $c_{j_0}$ is
elliptic by (\ref{5.5}) and (\ref{5.53}).
While the closed geodesic $c_\star$ is
elliptic by (\ref{5.4}) and (\ref{5.60}).

{\bf Case 2.} { \it If Case 3.2 holds.}

In this case, we can find two prime closed geodesics $c_{k_1}$ and $c_{k_2}$\
such that $i(c_{k_1})=n-1=i(c_{k_2})$. Then by the proof of
Theorem 5 in \cite{Rad4}, both $c_{k_1}$ and $c_{k_2}$ are
elliptic.

The proof of Theorem 1.3 is complete.

\medskip

\noindent {\bf Acknowledgements.} I would like to sincerely thank my
advisor, Professor Yiming Long, for introducing me to the theory of
closed geodesics and for his valuable help and encouragement
to me in all ways.

\bibliographystyle{abbrv}

\bigskip

\end{document}